\documentclass[a4paper,12pt,twoside]{article}

\usepackage[latin1]{inputenc}
\usepackage{amssymb}
\usepackage{amsmath}
\usepackage{amsthm}
\usepackage{amsfonts}
\input xy
\xyoption{all}

\newcommand{\remark}[1]{\begin{rem}\emph{#1}\end{rem}}
\newcommand{\Alb}{\hbox{Alb}}
\newcommand{\SAlb}{\hbox{SAlb}}

\newcommand{\bb}[1]{\mathbb{#1}}
\newcommand\calo{\mathcal O}
\newcommand{\mc}[1]{\mathcal{#1}}

\newcommand\CH{\operatorname{CH}}
\newcommand\h{\operatorname{h}}
\newcommand\HH{\operatorname{H}}
\newcommand\Br{\operatorname{Br}}
\newcommand\Gal{\operatorname{Gal}}

\newcommand\Pic{\operatorname{Pic}}
\newcommand\spec{\operatorname{spec}}

\input cyracc.def
\font\tencyr=wncyr10
\def\cyr{\tencyr\cyracc}

\theoremstyle{plain}
\newtheorem{X}{X}[section]

\newtheorem{T}{Theorem}
\newtheorem{TL}[X]{Lemma}

\newtheorem{D}{Definition}
\newtheorem{example}{Example}
\newtheorem{rem}[X]{Remark}

\newtheorem{Thm}[X]{Theorem}
\newtheorem{Prop}[X]{Proposition}

\newtheorem{Cor}[X]{Corollary}

\theoremstyle{definition}

\newcommand{\Proof}{\noindent\emph{Proof}\qquad }

\newcommand{\Z}{\mathbb{Z}}
\newcommand{\Q}{\mathbb{Q}}
\newcommand{\QZ}{\Q / \Z}
\newcommand{\F}{\mathbb{F}}
\newcommand{\Fp}{\F_p}

\newcommand{\sha}{\mbox{\cyr X}}
\newcommand{\ba}{\mbox{\cyr B}}
\newcommand{\bmo}{Brauer-Manin obstruction}

\newcommand{\NS}{\mathrm{NS}}

\newcommand{\ses}[5]{$$
\xymatrix@1{ 0 \ar[r] & {{#1}} \ar[r]^-{{#2}} & {{#3}} \ar[r]^-{{#4}} &
{{#5}} \ar[r] & 0
\\ }$$}

\newcommand{\sesdot}[5]{$$
\xymatrix@1{ 0 \ar[r] & {{#1}} \ar[r]^-{{#2}} & {{#3}} \ar[r]^-{{#4}} &
{{#5}} \ar[r] & 0. \\ }$$}

\newcommand{\sesbig}[5]{$$\xymatrix@1{{\raisebox{1.0ex}[3.0ex][1.0ex]{$0$}}
\ar@<0.6ex>[r] & {\raisebox{1.0ex}[3.0ex][1.0ex]{${#1}$}}
\ar@<0.6ex>[r]^{#2} & {\raisebox{1.0ex}[3.0ex][1.0ex]{${#3}$}}
\ar@<0.6ex>[r]^{#4} & {\raisebox{1.0ex}[3.0ex][1.0ex]{${#5}$}}
\ar@<0.6ex>[r] & {\raisebox{1.0ex}[3.0ex][1.0ex]{$0$}}}$$}

\newcommand{\from}{\colon}

\renewcommand{\implies}{\Longrightarrow}

\renewcommand{\TH}[3]{\widehat H^{#1}(#2, #3)}
\renewcommand{\H}[3]{H^{#1}(#2, #3)}

\newcommand{\Xbar}{\overline{X}}
\newcommand{\Xv}{X_v}
\newcommand{\Xvbar}{\overline{\Xv}}

\newcommand{\kv}{k_v}

\newcommand{\kvbar}{\overline{\kv}}

\newcommand{\Comment}[1]{}

\newcommand{\fixme}[1]{\textbf{#1}}
\begin{document}
\title{On the Brauer-Manin obstruction for zero-cycles on curves}
\author{Dennis Eriksson, Victor Scharaschkin}

\maketitle

\abstract{We wish to give a short elementary proof of S. Saito's
result that the Brauer-Manin obstruction for zero-cycles of degree
1 is the only one for curves, supposing the finiteness of the
Tate-Shafarevich-group $\sha^1(A)$ of the Jacobian variety. In
fact we show that we only need a conjecturally finite part of the
Brauer-group for this obstruction to be the only one. We also
comment on the situation in higher dimensions.} \tableofcontents

\newpage

Acknowledgements: The first author would like to thank D. Harari, who
was his advisor during the writing of this article, for inteteresting
remarks and discussions. He would also like to thank Colliot-Th{\'e}l{\`e}ne
and Olivier Wittenberg for stimulating conversations. He would
especially like to point out that Colliot-Th{\'e}l{\`e}ne suggested that an
early manuscript was a counter-example to the "global to
local"-principle; the overall idea was good but there were many local
obstructions. Both authors hope that the current manuscript is no more
in defiance of this principle. 

\section{Brauer-Manin obstruction}
We first recall the definition of the Brauer-Manin obstruction. Let
$X$ be any smooth variety over a field $k$. Set $\Br(X) :=
\HH^2_{et}(X, \bb G_m)$. By functortiality, an $L$-point (for $L/k$ a
finite field-extension) $\spec L \to X$ defines a homomorphism $\Br X
\to \Br L$. Furthermore, since $L/k$ is finite, we can take
corestriction $\Br L \stackrel{Cores}{\to} \Br k$. Hence, by extending
by linearity, we obtain a pairing of zero-cycles of $X$ and $\Br X$: 
$$
Z_0(X) \times \Br(X) \to \Br k.
$$
 Now, let $k$ be a number field, and set $k_v$ to be the
completion of $k$ at a place $v$. For a $k$-variety $X$, we denote
by $X_v = X \times_k k_v$, and by $\overline X = X \times_k
\overline k$ for a separable closure $\overline k$ of $k$. We define
the group of adelic zero-cycles $Z_0^a(X)$ as the subgroup of
$\prod_v Z_0(X_v)$ such that almost all zero-cycles $(z_v)_v$ extend
to zero-cycles over some model $\mc X \to \calo_v$ of $X_v \to k_v$
for $\calo_v$ the integer of $k_v$. If $X/k$ is proper, then $X$
admits a proper model over an open subset of $\spec \calo_k$, and
hence by base-change to our local fields and by the valuative
criterion of properness, these "rational" zero-cycles extend to
integral zero-cycles. If $X$ is
smooth, and has an adelic point, it is true that any open $U$ of $X$
also has an adelic point. This follows from the implicit function
 theorem, which is a consequence of the inverse function theorem
 \cite{Serre3}, Theorem 1, LG. 3.13. The same proof shows that if $X$ has an
adelic zero-cycle of degree 1, that is, adelic zero-cycles $(z_v)$
 where all $z_v$ have degree 1, then $U$ also does. \\
By virtue of $\Br \calo_w = 0$ (by local class field theory), the
following pairing is well-defined on the adelic zero-cycles
$Z_0^a(X):$
$$Z_0^a(X) \times \Br(X) \to \bb Q/\bb Z,$$ via
$$((z_v), A) \mapsto \langle (z_v), A \rangle = \sum_v i_v(A(z_v)),$$ where $i_v \from
\Br k_v \to \bb Q/\bb Z$ is the usual injection (which is an
isomorphism for finite places) given by local
class field theory. Denote by $\bb A_k$ the $k-$adeles. A necessary
condition for the non-emptieness of $X(k)$ is obviously the
non-emptieness of $X(\bb A_k)$, and if there is a class of varieties
which satisfy the converse of this statement, we say they satisfy the
Hasse principle. It is however wellknown that this class does not
contain all varieties, a famous counter-example due to Selmer is given by the
curve of genus 1 given by
$$3x^3 + 4y^3 + 5z^3 =0.$$
Now, define the set
$$X(\bb A_k)^{\Br X} = \left \{ (z_v) \in X(\bb A_k)\mid \langle
(z_v), A \rangle = 0, \forall A \in \Br X \right \}.$$
Manin~\cite{Manin} noted that the set of rational points is included
in this set via the global reciprocity law of class field theory,
and suggested this as a stronger Hasse condition. We define, for $B
\subseteq \Br X$ (or $\Br X/\Br k$, the pairing factors over the
image of $\Br k \to \Br X$. Here and henceforth $\Br X/ \Br k$ means the quotient by
the image of $\Br k$.)
$$Z^a_0 (X)^B = \left \{(z_v) \in Z^a_0 (X) \mid
\forall v, \deg z_v = 1, \langle (z_v), A \rangle = 0, \forall A \in B ,
\right \}.$$
Using that the corestriction-map is the identity on the level of
$\bb Q/\bb Z$ for non-archmidean places (see \cite{Serre}, XI, Prop 2,
ii) and XIII, Theorem 1) and the
fundamental short exact sequence of class-field theory
$$0 \to \Br k \to \oplus_v \Br k_v \to \bb Q/\bb Z \to 0,$$
one shows that, $Z_0(X)^1$, the zero-cycles of degree 1, are indeed
included in this set. If this set is non-empty, we say there is no
obstruction associated to $B$ for existence of zero-cycles of degree
1. Now let $\mc V$ be a class of varieties. If for all $X$ in $\mc
V$ we have
$$Z_0^a(X)^{B} \neq\emptyset \implies Z_0(X)^1 \neq 0$$ then we
say the Brauer-Manin obstruction is the only one to the existence of
zero-cycles of degree 1 associated to $B$. 
Let $\Br_1(X) = \ker [\Br X \to \Br (\overline X)].$ The group
$$\ba(X) :=  \ker \left[\Br_1 X / \Br k \to \prod_v \Br_1 X_v / \Br k_v \right]$$ is the
group of (algebraic) locally constant elements of the Brauer group,
which is canonically isomorphic to $$ \sha^1(\Pic(\overline X)) =
\ker \Big(\HH^1(k, \Pic(\overline X)) \to \prod_v \HH^1(k_v,
\Pic(\overline X_v))\Big )$$ via the Hochschild-Serre spectral
sequence (see \cite{torsors}, this is an isomorphism essentially
because $H^3(k, \bb G_m) = 0$ for local and global fields). Here and
henceforth $H^i(k, M)$ denotes {\'e}tale cohomology, which reduces
to Galois-cohomology of $\Gamma_k = \Gal(\overline k/k)$ with
values in $M$. \\

Let $X$ be a smooth quasi-projective variety defined over a number
field $k$. If $X$ is also proper, denote the Albanese variety by
$\Alb^0_X$ or $A$ and the Picard variety by $B$.  The \emph{index}
$I = I_X$ of a variety $X$ is defined to be the cardinality of the
cokernel of (see the next section for definition of Suslin homology)
$\deg \from \h_0(X) \to \bb Z$ or $\deg: Z_0(X) \to \bb Z$. The
\emph{period} $P = P_X$ of a variety is defined as the cardinality
of the cokernel of $\deg \from \h_0(\overline X)^{\Gamma_k} \to \bb
Z$. Define the generic period $\widetilde P$ as the supremum of all
$P_U$ for all open subsets $U$ of $X$. Note that $P_U \mid I_U$ and
and it is wellknown that the index of an open subset of $X$ is the
same as that of $X$, we see that all $P_U$ are bounded by $I$, so
the supremum exists. Moreover $\widetilde P | I$.

Colliot-Th{\'e}l{\`e}ne~\cite{CT} has conjectured that the \bmo\ is
the only obstruction to the existence of global $0$-cycles of degree
1 on $X$. We shall prove a very weak version of this conjecture.

\begin{Thm}\label{main}
Assume that $$\sha^1(A) = \ker[\HH^1(k,A) \to \bigoplus_v
\HH^1(k_v,A)]$$ is finite. If $Z^a_0(X)^{\ba(X)} \neq \emptyset$
then $\widetilde P=1$.
\end{Thm}
Thus the \bmo\ is the only obstruction to the global generic
\emph{period} being 1. A related result is the following stronger
version of a theorem originally due to Saito~\cite{Saito} (stronger,
because here we only need the conjecturally finite group $\ba (C)$
as opposed to the whole Brauer group). Saito's theorem has also been
reproven by Colliot-Th{\'e}l{\`e}ne in  \cite{CT}.

\begin{Thm}\label{SaitoThm} Let $C$ be a smooth projective curve
over a number field $k$, let $A$ be its Jacobian, and assume that
$\sha^1(A)$ is finite. Then the obstruction associated to $\ba(C)$
for zero-cycles of degree 1 is the only one:
$$\hbox{if } Z_0^a(X)^{\ba(C)} \neq\emptyset \hbox{ then } I =
1.$$
\end{Thm}

We need a couple of lemmas and some notation.  The
\emph{semi-Albanese torsor} of $X$ will be a morphism $p \from X \to
\SAlb^1_X$, where $\SAlb^1_X$ is a torsor under the semi-Albanese
variety $\SAlb^0_X$, with a universal property with respect to
torsors under semi-abelian varieties. For the ones only interested
in \ref{SaitoThm}, we urge you to skip the next section and jump
immediately to section 3.

\section{Construction of the semi-Albanese Torsor} Let $U$ be a
quasi-projective, smooth variety over a perfect field $k$. Recall
that a semi-abelian variety is a commutative group-variety which is
an extension of an abelian variety by a torus. Suppose we are given
a morphism $p: U \to S$ where $S$ is a torsor under a semi-abelian
variety, with the following universal property: Given any morphism $f: U \to T$, where $T$ is a torsor under a semi-abelian
variety, there is a unique morphism $ S \to T $ such that $h p = f$.
This clearly determines the tripple $(U, S, p)$ up to unique
isomorphism, and we will call it the semi-Albanese torsor of $U$.
The principal result of this section is Theorem \ref{period} which related the
period of $U$ to an invariant depending only on the semi-Albanese
torsor. \\

The following is a formal consequence of its solution over an
algebraically closed field 
\cite{Serre1} (see Theorem 7) and the descent theory of \cite{Serre2} (see p. 112, 4.22).

\begin{T} Let $U$ be a quasi-projective smooth variety over a perfect field
$k$, then the semi-Albanese torsor exists.
\end{T}

\begin{proof} If we are over an algebraically closed field, the
theorem is Theorem 7 of \cite{Serre1}. Let $\overline{k}$ be an
algebraic closure of $k$, and denote by $\SAlb_{\overline{U}}$ the
semi-Albanese variety which solves the universal problem for
$\overline U$ over $\overline k$. Also denote by, for a $\overline
k$-variety $V$, $V^\sigma$ the base-change of $V \to \overline k$ and
$\overline k \stackrel{\sigma}{\to} \overline k$. If $\sigma \in \Gamma_k$, if we let $ p^\sigma =
\sigma(p(\sigma^{-1}()) $, there exists by universality a unique $\overline
k$-morphism $h_\sigma$ such that
$$\xymatrix{ \overline U \ar[dr]^{p^\sigma} \ar[r]^p & \SAlb_{\overline{U}} \ar[d]^{h_\sigma} \\
& \SAlb_{\overline{U}}^\sigma}$$ commutes. Because it is unique, it
automatically satisfies the cocycle-condition $$h_{\tau \sigma} =
h_{\tau}^{\sigma} h_{\sigma}.$$ By the theory of
descent \cite{Serre}, p. 108, 4.20. this descends the variety to a
variety $\SAlb_U^1$ defined over $k$, and also descends the
$\overline k$-morphism $p$ to a $k$-morphism $p: U \to \SAlb_U^1$. The morphism $h_\sigma$ moreover factors uniquely
as a translation followed by an isomorphism of group-schemes. If we
denote by $h^0_\sigma$ this group-scheme isomorphism, it also
automatically satisfies a cocycle condition. If we descend $\SAlb_{\overline{U}}$
by this cocycle it also descends to a $k$-variety $\SAlb_U^0$, and
$\SAlb^1_U$ is easily verified to be a torsor under $\SAlb_U^0$.

Now, to verify the universal property, suppose that we are given a
$k$-morphism $f$ from $U$ to a torsor under a semi-abelian variety
$T$. Over the algebraic closure, we choose a point in $T$ so it
becomes a semi-abelian variety, which we also denote by $\overline
T$, and via the isomorphism $\SAlb_{\overline{U}} \simeq
\SAlb_{\overline{U}}^1$ we obtain a unique $\overline k$-morphism
$h: \SAlb_{\overline{U}}^1 \to \overline T$, such that $h p = f$.
However, if $\sigma \in \Gamma_k$, then $hp = f = f^\sigma =
h^\sigma p^\sigma = h^\sigma p$, and since $h$ was unique $h =
h^\sigma$, so $h$ descends to a $k$-morphism $U \to T$ such that $f
= hp$.

\end{proof}

Recall the following construction of Suslin homology: An elementary
finite correspondance from $Y$ to $X$ over a field $k$ is an
integral subscheme $Z$ of $X \times_k Y$ which is finite and
surjective over $Y$. A finite correspondance between $X$ and $Y$ is
a formal $\bb Z$-linear sum of elementary finite correspondances,
and we denote the group of such as $Cor(Y, X)$. Note that the finite
correspondances from $k$ to $X$ is just the group of zero-cycles on
$X$.

Denote by $i_0$ and $i_1$ the points $0$ and $1$ in $\bb A^1$. Given
a finite correspondance from $\bb A^1$ to a variety $X$, we get an
associated finite correspondance from $k$ to $X$  via
$$Z \mapsto i^*_0 Z - i^*_1 Z. $$
We define $\h_0(X)$, the 0-th Suslin homology of $X$, to be the
group of zero-cycles on $X$ modulo the group generated by finite
correspondances coming from $\bb A^1$ to $X$ in the above sense. We
note the following properties, which are not difficult to show:
\begin{enumerate}
    \item If the structural morphism $X \to \spec k$ is proper, then $\h_0(X) =
    \CH_0(X)$, the $0$-th Chow group.
    \item $h_0$ is covariant with respect to morphisms between
    varieties.
    \item The structural morphism $X \to \spec k$ induces the degree
    map $\deg: \h_0(X) \to \h_0(k) = \bb Z.$
\end{enumerate}

\remark{If $X$ is also proper, then the semi-Albanese variety is the
Albanese variety, and the semi-Albanese torsor is an "Albanese
torsor", i.e. it is universal with respect to morphisms into torsors
under abelian varieties. In this case we write it as $\Alb_X^1$
instead. If $H^1(X, \calo_X) = 0$ the abelian-variety part of the
semi-Albanese variety is trivial and the semi-Albanese torsor is a
torsor under a torus and is universal with respect to morphisms to
torsors under tori. }

The semi-Albanese scheme of $X \to \spec k$ is the $k$-group scheme
$$\underline{\SAlb}_{X/k} = \coprod_{n \in \bb Z} \SAlb^n_X$$ where
$\SAlb^n_X$ is the $n$-fold Baer sum of torsors, and for $n = 0$ it
is the semi-Albanese variety. In \cite{duality},  1.2, Ramachandran
shows this is a group-scheme with functorial and universal
properties, which we choose to reproduce somewhat. We have an
obvious $\Gamma_k$-equivariant map $X(\overline{k}) \to
\SAlb^1_X(\overline{k})$, which gives in the natural way a map from
the group of zero-cycles $Z_0(\overline {X})$ to the $\overline
k$-points of the Albanese scheme. Taking Galois-invariants gives a
group-homomorphism from $Z_0(X) \to \underline{\SAlb}_{X/k}(k)$. By
\cite{Albanese}, Lemma 3.1, this factors over the group $h_0(X)$.

All in all, there is a canonical homomorphism $\h_0(X) \to
\underline{\SAlb}_{X/k}(k)$ such that the restriction to degree 0 is
the generalized Albanese map  of $\cite{Albanese}$ and the
structural morphism $X \to \spec k$ induces the following
commutative diagram
$$\xymatrix{0 \ar[r] & A_0(X) \ar[d] \ar[r] & \h_0(X)
\ar[r]
\ar[d] & h_0(k) = \bb Z \ar@{=}[d]  \\
0 \ar[r] & \SAlb_X^0(k) \ar[r] & \ar[r] \underline{\SAlb}_{X/k}(k) &
\underline{\SAlb}_{k/k}(k) = \bb Z }$$

Now, any torsor $T$ under an semi-abelian variety $A$ defines an
element in $\HH^1(k,A)$, which we denote by $[T]$. This element is
trivial exactly when $T$ has a $k$-point, and is isomorphic to $A$
over $\overline{k}$. Define the period of $T$ to be the order of the
element $[T]$. This is compatible with our previous definition:

\begin{Thm}\label{period} Let $X$ be a quasi-projective, smooth variety over a
perfect field $k$. Then the period of $X$ is equal to the period of
the semi-Albanese torsor.
\end{Thm}

\begin{proof} By the preceding remarks we have
the following commutative diagram of Galois modules with exact rows:
$$\xymatrix{0 \ar[r] & \ar[r] \ar[d]^p A_0(\overline X)
& \ar[r] \ar[d] \h_0(\overline X) & \ar[r] \ar@{=}[d] \bb Z & 0 \\
0 \ar[r] & \ar[r] \SAlb_X^0(\overline k) & \ar[r]
\underline{\SAlb}_{X/k}(\overline k) & \bb Z \ar[r] & 0}.$$
Here $p$ is the generalized Albanese map (loc. cit.). Taking Galois cohomology
gives us the following diagram
$$\xymatrix{ \h_0(\overline X)^{\Gamma_k} \ar[d] \ar[r] &
\bb Z \ar@{=}[d] \ar[r] & \HH^1(k, A_0(\overline X)) \ar[d]^p \\
\underline{\SAlb}_{X/k}(k) \ar[r] & \bb Z \ar[r] & \HH^1(k,
\SAlb_X^0(\overline k))}$$
The image of 1 in $\HH^1(k, \SAlb^0_X)$ is represented by the
cocycle $\sigma \mapsto x_0^\sigma - x_0$, for $x_0 \in
\SAlb^1_X(\overline k)$, i.e. the class of the Albanese torsor. Now,
the generalized theorem of Roitman says that the kernel of the
Albanese map is uniquely divisible (see main theorem of \cite{Albanese}), and
$\HH^i(k, \bb Q) = 0$ for $i > 0$, so that the rightmost
homomorphism is an isomorphism. A diagram-chase finishes the lemma.
\end{proof}

\remark{Notice also that for the case of a proper curve the lemma is
trivial. Indeed, the (semi-)Albanese scheme is the Picard scheme,
and all the maps involved are canonical isomorphisms; there is
nothing to prove. Also, in this case $\SAlb^1_X$ is usually written
as $\Pic^1_X$.}

One might ask what kind of relationship the index and period has
in general for torsors under abelian varieties. In this direction
we have the following general proposition which is well known:
\begin{Prop} \label{periodforab} Let $T$ be a torsor under an abelian variety $A$ of dimension $d$ over a
number field $k$. Then
$$P \mid I \mid P^{2d}.$$
\end{Prop}
\begin{proof} The first divisibility $P \mid I$ is classical, and not difficult. Indeed, it is
clear that if $T$ has a $k$-point, then $P = I = 1$. The general
case follows from a restriction-corestriction-argument with respect
to the zero-cycle which gives the index. For the second, we have the
diagonal map or multiplication by $P$-map $T \to^P T^P \simeq A$,
where $T^P$ is the $P$-fold Baer-sum of torsors. This is of degree
$P^{2d}$, and so we obtain a zero-cycle of degree $P^{2d}$ in $T$
above a rational point in $A$.
\end{proof}

\section{Proofs of the Main Theorems}

In this section we give the proofs of the two main theorems,
roughly as follows. First we show that under the right conditions the
period is equal to 1 (Lemma \ref{torsab}), and we show that for a
curve that local conditions give that the period is actually equal to
the index (Lemma \ref{curve}). We then prove that the locally constant
elements is invariant under restriction to Zariski-open (Lemma
\ref{Lemma:locallyconstant}) and then put all of this together to
prove the theorems.

\begin{TL} \label{torsab}Let $V$ be a torsor under an semi-abelian variety
$S$ which is an extension of an abelian variety $A$ by a torus $T$,
defined over a global field $k$ and suppose that $\sha^1(A)$ is
finite. Then the obstruction associated to $\ba(V)$ for zero-cycles
of degree one is the only one for points on $V$. That is, if
$Z_0^a(V)^{\ba(V)} \neq\emptyset$ then $V(k) \neq \emptyset$.
\end{TL}

\begin{proof} The statement is well known if we replace zero-cycles with
points and take $S$ to be an abelian variety (see Theorem 5.2.3 of \cite{torsors} for
a proof, or the original article of Manin \cite{Manin}) or $S$ to be a
torus (see Theorem 5.2.1, \cite{torsors}). Since $V$
has a zero-cycle of degree one locally everywhere, it has local
period 1, so actually has a $k_v$-point by Theorem~\ref{period} or the
first part of the argument in the proof of Proposition
\ref{periodforab}. For each place $v$, let $Q_v$ be such a point, and
suppose that $Q = (Q_v)_v$ is adelic. A
restriction--corestriction argument shows that $i_v(\alpha(z_v)) =
\deg(z_v)i_v(\alpha(Q_v))$ for any zero-cycle $z_v$ on $X_v$, for
$\alpha$ locally constant. Hence we can replace all zero-cycles of
degree one with local points, and this adelic point $Q$ will be orthogonal to
$\ba(V)$. The result for $\ba(V)$ in the case of $S$ neither a torus
or an abelian variety is an
unpublished result of Harari and Szamuely (\cite{Harari}). 
\end{proof}

\begin{TL}\label{curve} Let $X$ be a smooth, proper curve over a global field $k$ and
assume that $X$ has a zero-cycle of degree one locally everywhere.
Then the index is equal to the period, $I = P$.
\end{TL}
\begin{proof} This proof can be found in Prop 2.5 of~\cite{MilneBTS}.
We include it for completeness. Writing out the lower terms of the
Hochschild-Serre spectral sequence
$$
E^{p,q}_2 = \HH^p(k,\HH^q_{et}(\overline{X}, \bb G_m) ) \Longrightarrow
\HH^{p+q}_{et}(X, \bb G_m)
$$
gives us the exact sequence
$$0 \to \Pic(X) \to \HH^0(k,\Pic(\Xbar)) \to \Br k \stackrel{j}\to \Br X.$$
If $X$ has a $k$-point, this point splits the map $j$ and so $j$ is
injective. By a restriction--corestriction argument the same stays true if
$X$ has a zero-cycle of degree 1. Global class-field theory tells us that
the map $\Br k \to \oplus_v \Br k_v$ is injective (the Hasse principle for
Severi-Brauer varieties). The condition that we have a zero-cycle locally
everywhere gives us that $\oplus_v \Br k_v \to \oplus_v \Br X_v$ is
injective, and one deduces that $\Pic(X) \to H^0(k, \Pic(\overline X))$
must be surjective. Then from the definition of index and period it is
obvious that they must be equal.
\end{proof}

The following corollary may be useful in calculations with curves.
\begin{Cor} Let $X$ be a smooth, proper curve over a global field $k$ of genus
$g \geq 2$ with points everywhere locally and Jacobian $A$. If
$\sha^1(A)[p]=0$ for all primes dividing $2g-2$ then the curve has
a rational $0$-cycle of degree $1$.
\end{Cor}
\Proof By lemma~\ref{curve} and lemma ~\ref{torsab}, the index of
the curve is equal to the order of the Albanese torsor in
$\sha^1(A)$. The existence of the canonical divisor implies $I \mid
(2g-2)$, and the hypothesis thus implies $I=1$.\qed

For example, if $g=2$ or $3$ and $\sha^1(A)[2]=0$ then $X$ has a
$0$-cycle of degree 1.  See~\cite{Flynn} for a discussion of how to
find curves of low genus with index 1 but no global points. \\

\begin{TL}\label{Lemma:locallyconstant} Suppose $U$ is open in $X$, then $$\ba(U) = \ba(X).$$
\end{TL}

\begin{proof} It is well-known that $\Br X$ stays the same if one removes a subvariety of codimension 2 or
more. Let $Z$ be a codimension one variety
of $X$, $U = X \setminus Z$. We want to show that $A \in \ba(U)$
is unramified at $Z$, i.e. that the residue $d(A)$ in $H^1(k(Z),
\bb Q/\bb Z)$ is 0, since then it extends to a class in $\ba(X)$.
First of all, since $A$ is everywhere locally constant, it is
everywhere locally unramified, so $A$ goes to 0 in $H^1(k(Z)
\times_k k_v, \bb Q/\bb Z)$. Choose an embedding of $k_v$ into
$\overline {k}$, and thus a map from $H^1(k(Z) \times_k k_v, \bb
Q/\bb Z) \to H^1(k(Z) \times_k \overline {k}), \bb Q/\bb Z)$ and
so the map $H^1(k(Z), \bb Q/\bb Z) \to H^1(k(Z) \times_k
\overline{k}, \bb Q/\bb Z)$ maps $d(A)$ to zero. Let $L$ be the
algebraic closure of $k$ in $k(Z)$. Then the above means that the
residue $d(A)$ comes from $H^1(L, \bb Q/\bb Z)$. The maps $H^1(L,
\bb Q/\bb Z) \to H^1(k(Z), \bb Q/\bb Z)$ and $H^1(L \times_k k_v,
\bb Q/\bb Z) \to H^1(k(Z) \times_k k_v, \bb Q/\bb Z)$ are
injective by the Hochschild-Serre spectral sequence (loc.cit.), so
$d(A)$ belongs to
$$\sha^1(L, \bb Q/\bb Z):= \ker \left[H^1(L, \bb Q/\bb Z) \to
\prod_v H^1 (L_v, \bb Q/\bb Z)\right ].$$ However, this group is
the direct limit of $\sha^1(L, \bb Z/n)$ and this group is zero.
Indeed, it is the group of cyclic extensions of order $m$ dividing
$n$ which are split everywhere. However, by the Chebotarev density
theorem, there are no such extensions, since the set of primes
which split completely have analytic density $1/m$. Hence $A$ is
unramified at $Z$ and extends to an element of $\ba(X)$.
\end{proof}

\begin{proof} (of Theorem~\ref{main} and \ref{SaitoThm}). Suppose
that $X$ is proper and smooth (we can do so by taking a smooth
compactification, by Hironaka say). Let $U$ be an open of $X$, and
let $p \from U \to \SAlb^1_U$ be its semi-Albanese torsor. Since
$\ba(U) = \ba(X)$ and these elements are locally constant, $U$ has
no $\ba(U)$-obstruction. By the projection-formula, the same holds
true for $\SAlb^1_U$. Because of the finiteness assumption of
$\sha^1(A)$, Lemma ~\ref{torsab} implies that the torsor is trivial
and so Lemma~\ref{period} says that its period is equal to one. This
is true for any open in $X$, and hence the generic period is 1. For
Theorem \ref{SaitoThm}, notice that in this case Lemma ~\ref{curve}
says that $I = P$, so $I = P = 1$. This proves the two main theorems. 
\end{proof}

%------------------------------------------------------------------------
%------------------------------------------------------------------------
%------------------------------------------------------------------------
%------------------------------------------------------------------------

\begin{Cor} Let $X$ be a smooth quasi-projective $k$-variety (of
arbitrary dimension larger than or equal to 1). If there is any open
$U$ in $X$ such that $\h_0(U) \to \h_0(\overline U)^{\Gamma_k}$ is
surjective and $\sha^1(A)$ is finite, then the \bmo\ is the only
obstruction to the existence of global $0$-cycles of degree 1 on
$X$.
\end{Cor}

Indeed, the surjectivity of this map implies that the index is equal
to the period. \\

\remark{We note that $P_U$ can indeed be larger than $P_X$ for $U$
open in $X$. For example, if $X$ is a proper curve of genus 0, then
via the anti-canonical embedding it can be written as a conic in
$\bb P^2$: $$X: aX^2 + bY^2 = cZ^2.$$ Hence the index is either 1 or
2, and it is 1 exactly when we have a rational point. Now, removing
two points at infinity, we obtain $$U: ax^2 + by^2 = c$$ which is a
torsor under a torus. Because $P_U$ divides $I$, it is either 1 or
2, and because the torsor is trivial exactly when $P_U$ is 1, we see
that $P_U = I$. Hence we have in this case that $\widetilde P = I$.
However since the Albanese of $X$ is trivial, $P_X$ is certainly
1. The same argument for any compactification of torsors under
tori allows us to recover a result by Colliot-Th{\'e}l{\`e}ne and Sansuc
saying that the Manin-obstruction is the only one for
smooth compactifications of $k$-torsors under tori (see \cite{torsors},
Theorem 5.3.1, and the remark afterwards saying that we only need to
consider locally constant elements). 
In any case, the generic period contains more information than the period
associated to only $X$. An interesting question  (suggested by
Colliot-Th{\'e}l{\`e}ne) would be to calculate
the generic period of (a compactification of) a non-abelian algebraic group and compare it to
its index. }

Suppose henceforth that $X$ is proper. In section~\ref{extra} we use
Tate duality to give a more explicit description of the set
$Z^a_0(X)^{\Br(X)}$ in some situations.

Let $\langle, \rangle_{\mathrm{Tate}} \from \sha^1(A) \times
\sha^1(B) \to \QZ$ denote the Cassels-Tate pairing. We now show
that the Brauer-Manin pairing can be related to this pairing using
the Albanese torsor.

>From Hochschild-Serre we have:
\begin{equation}
\label{snake1} \xymatrix{& & & 0 \ar[d] \\ & 0 \ar[d] & & \sha^1(\Pic(\overline X)) \ar[d] \\
0 \ar[r] & \Br(k) \ar[r]^{i}\ar[d] & \Br_1(X) \ar^{r}[r]\ar^{s}[d] &
\H{1}{k}{\Pic(\Xbar)} \ar[d]
\ar[r] & 0 \\
0 \ar[r] & \prod_v \Br(\kv) \ar[r]^{i_v} & \prod_v \Br_1(\Xv) \ar[r] &
\prod_v \H{1}{\kv}{\Pic(\Xvbar)} \ar[r] & 0 \\ }
\end{equation}
Let
$$\phi \from \sha^1(\Pic(\overline X)) \to
\QZ$$ be the map provided by the snake lemma. Note that $\phi$ is
defined using (any) sections to the $i_v$ (which is why it maps
into the cokernel of $\Br(k) \to \bigoplus \Br(\kv)$). Thus for
any local $0$-cycles $(c_v)$ we may calculate $\phi(b)$ using the
formula $\phi(b) = \langle (c_v), r^{-1}(b)
\rangle_{\mathrm{BM}}$. Here $r^{-1}(b)$ denotes the inverse image
in $\Br_1(X)/\Br k$, and we know that constant elements don't
cause obstruction, so the pairing makes sense.

Let $\rho$ be the natural map $\rho \from \sha^1(B) \to
\sha^1(\Pic(\overline X))$. Let
$$
T =r^{-1} \left(\rho(\sha^1(B)\right)) \subseteq \Br_1(X).
$$
Then for any $b$ in $T$, $\langle (c_v), b \rangle_{\mathrm{BM}} =
\phi(r(b))$ (independent of the choice of local $0$-cycles). Since $\phi
\circ r = 0$ iff $\phi = 0$, we have

\begin{equation}
\label{relatetophi} Z^a_0(X)^{T} \neq \emptyset \iff \phi = 0.
\end{equation}

\begin{Thm}\label{CTpair} Let $[X] \in \sha^1(A)$ represent the
Albanese torsor, $p: X \to \Alb_X^1$ the Albanese torsor morphism.
Then for all $b \in \sha^1(B)$,
\begin{equation}
\langle (c_v), p^* b \rangle_{\mathrm{BM}} = \langle p_*(c_v), b
\rangle_{\mathrm{BM}}  = \langle [X],b \rangle_{\mathrm{Tate}},
\end{equation}
\end{Thm}
\Proof This is similar to Manin's theorem for genus 1 curves and
follows from the homogeneous space definition of the Cassels-Tate
pairing~\cite[Poonen \& Stoll section 3.1]{PoonenStoll}
and~\cite[Milne, I.6.11]{MilneAD}, or directly using the
definition of the map $X \to \Alb_X^1$ and the projection
formula.\qed

%---------------------------------------------------------------------
%----------------------------------------------------------------------
\section{A description of the Brauer Set}
%----------------------------------------------------------------------
%----------------------------------------------------------------------

\label{extra} In this section we suppose $X$ proper and we give a
description amenable to computation of the Brauer-Manin set
$Z^a_0(X)^{\Br_2(X)}$ (and related objects) for a certain subgroup
$\Br_2(X) \subseteq \Br(X)$. Let $B$ be the Picard variety of $X$,
and let $\widehat{H}^0$ denote Tate cohomology (see \cite{Serre}).
Assume from now on that the map $\H{0}{k}{\Pic(\Xbar)} \to
\H{0}{k}{\NS(\Xbar)}$ is surjective, so that so that we may regard
$\H{1}{k}{B}$ as a subgroup of $\H{1}{k}{\Pic(\Xbar)}$. Let
$\Br_2(X) \subseteq \Br_1(X)$ be the preimage of this subgroup under
the map $r \from \Br_1(X) \to \H{1}{k}{\Pic(\Xbar)}$ arising in the
Hochschild-Serre spectral sequence.  Thus $r$ induces an isomorphism
$$
\overline{r} \from \frac{\Br_2(X)}{\Br(k)} \simeq \H{1}{k}{B}.
$$

Assume that $X$ has local period $1$ at every place and that
$\sha^1(A)$ is finite.  There are local Albanese maps $j_v \from
\TH{0}{\kv}{A_0(\Xvbar)} \to \TH{0}{\kv}{A_v}$. Define for an
abelian group $G$, $G^* = \hom(G, \bb Q/\bb Z)$, and $\widehat G$
the profinite completion of $G$. Then from the Tate duality
sequence for $A$ (see \cite{MilneAD}, Proposition 6.23(b) and
Remark 6.14(a) ) we can form the following diagram:

\begin{equation}\label{ExplicitDiagram}
\xymatrix{ 0 \ar[r]& C
 \ar[r]\ar[d] &\prod_v \H{0}{\kv}{A_0(\Xvbar)}
\ar[r]^-{\beta}\ar[d]_j & \left(\frac{\Br_2(X)}{\Br(k)}
\right)^{*}\ar[d]_{\wr} \\ % \wr is wreath product sign to give vertical tilde isomorphism symbol
0 \ar[r] &\widehat{A(k)} \ar[r] &\prod_v \TH{0}{\kv}{A_v} \ar[r] &
\H{1}{k}{B}^{*} }
\end{equation}
where $j=(j_v)$, $\beta$ is the map making the diagram commute
and $C = \ker \beta$.   Thus
\begin{equation}\label{ImageOfC}
j(C) = \widehat{A(k)} \cap j(\H{0}{\kv}{A_0(\Xvbar).}
\end{equation}
Here $C$ can be viewed as a Brauer-Manin type set of $0$-cycle
classes as detailed below.

Suppose $\overline{z}$ is a 0-cycle rationally equivalent to $0$
over $\kvbar$. Then there exists a finite extension $L_w$ of $k_v$
such that $\overline{z_v}$ is defined over $L_w$ and the rational
equivalence is also defined over $L_w$.  Thus $\overline{z}$ pairs
to 0 in $\Br L_w$ under the pairing $Z_0(X_w) \times \Br (X_w) \to
\Br L_w$.   Hence we can define a pairing on $0$-cycle classes
$\H{0}{\kv}{\CH_0(\Xvbar)} \times \Br (\Xv) \to \Q/\Z$. For $T
\subseteq \Br(X)$ let
$$
C_0^a(X)^T = \{ [z_v] \in \prod \H{0}{\kv}{\CH_0(\Xvbar)} \mid
\deg(z_v) = 1, \quad \langle z_v, b \rangle = 0 \quad \forall b
\in T\}.
$$
This is the Brauer-Manin obstruction to the global period being
$1$.  If the period is 1, by translating we may identify
$0$-cycles classes of degree $1$ with classes of degree $0$. Then
$C_0^a(X)^{\Br_2(X)}$ gives rise to the set $C$ inside $\prod_v
\H{0}{\kv}{A_0(\Xvbar)}$.  Thus (\ref{ImageOfC}) gives a
description, not quite of $C_0^a(X)^{\Br_2(X)}$ but its image
after mapping into the $A_v$.

If $I=1$ then we may identify $0$-cycles of degree $1$ with those
of degree 1 and replace $\H{0}{\kv}{ A_0(\Xvbar)}$ by $Z_0(X_v)^1$
in diagram~\ref{ExplicitDiagram}.  (Abusing notation, w e continue
to call the map into the Albanese $j$). Similarly, if $I=1$ we may
translate $k$-points to $0$-cycles of degree 0.

In summary:

\begin{Thm}
\label{toomanyhypotheses} Suppose $P_v=1$ for all $v$, that the
map $\lambda \from \H{0}{k}{\Pic(\Xbar)} \to \H{0}{k}{\NS(\Xbar)}$
is surjective and $\sha^1(A)$ is finite.
\begin{enumerate}

\item Suppose $P=1$.  Let $C$ be the image of
$C_0^a(X)^{\Br_2(X)}$ inside $\prod_v \H{0}{\kv}{A_0(\Xvbar)}$.
Then $j(C) = \widehat{A(k)} \cap j(\H{0}{\kv}{A_0(\Xvbar)}$ inside
$\prod \TH{0}{\kv}{A_v}$.

\item If $I=1$ then we have similarly, after translation by a
$0$-cycle class of degree 1, $j(Z_0^a(X)^{\Br_2(X)}) =
\widehat{A(k)} \cap j(\prod Z_0(\Xv)^1)$ inside $\prod
\TH{0}{\kv}{A_v}$ and $j(X^{\Br_2(X)}(\bb A_k)) = \widehat{A(k)}
\cap j(X(\mathbb{A}_k))$ inside $\prod \TH{0}{\kv}{A_v}$.

\item{If $\H{1}{k}{\NS(\Xbar)} = 0$ then $\Br_2$ may be replaced
with $\Br_1$.}
\end{enumerate}
\end{Thm}

This description is probably most easily applied to prove that
certain Brauer-Manin sets are empty, by testing if the
intersection is empty.  If $j$ is injective we obtain a complete
description of the part of the Brauer-Manin set coming from
$\Br_2(X)$.

Note that curves automatically satisfy the condition on $\lambda$,
and $\Br_2(X) = \Br(X)$.  Moreover $j$ is injective so we obtain a
complete description of $X(\bb A_k)^{\Br X}$ in this case.

%---------------------------------------------------------------------------
%
% *** Bibliography
%
%---------------------------------------------------------------------------
\pagebreak 
\end{document}